\documentclass[7pt]{article}
\usepackage{float}
\usepackage{graphicx,epstopdf}
\usepackage{bbm}
\usepackage{amsfonts}
\usepackage{lipsum}
\usepackage{mathrsfs}
\usepackage{amsmath}
\usepackage{csquotes}
\usepackage{dsfont}
\usepackage{amssymb}
\usepackage{theorem}
\usepackage{graphicx}
\usepackage[dvips]{epsfig}
\usepackage{latexsym}
\usepackage{exscale}
\usepackage[latin1]{inputenc}

\newtheorem{thm}{Theorem}[section]
\newtheorem{lem}{Lemma}[section]
\newtheorem{proof}{Proof}[section]
\newtheorem{prop}{Proposition}[section]

\newtheorem{defn}{Definition}[section]

\newcommand\blfootnote[1]{%
  \begingroup
  \renewcommand\thefootnote{}\footnote{#1}%
  \addtocounter{footnote}{-1}%
  \endgroup
}
\numberwithin{equation}{section}
\begin{document}
\small
\title{Representation Theorems of $\mathbb R$-trees and Brownian Motions Indexed by $\mathbb R$-trees}
\date{\vspace{-5ex}}

\maketitle
\begin{center}
\small{Asuman G\"{u}ven AKSOY, Monairah AL-ANSARI and Qidi PENG}
\end{center}
{\bf Abstract.} We provide a new representation of an $\mathbb R$-tree by using a special set of metric rays. We have captured the four-point condition from these metric rays and shown an equivalence between the $\mathbb R$-trees with radial and river metrics, and these sets of metric rays. In stochastic analysis, these graphical representation theorems are of particular interest in identifying Brownian motions indexed by $\mathbb R$-trees.

\section{Introduction}
One of the central object in stochastic analysis is Brownian motion,  which is the microscopic picture emerging from a particle moving in $n$-dimensional space and the nature of Brownian paths is of special interest. For example, the Brownian motion $B$ indexed by Euclidean space $(\mathbb R,|\cdot|)$ has stationary independent increments, i.e., $B(x_2)-B(x_1)$ and $B(x_4)-B(x_3)$ are independent and equally distributed if $x_1<x_2\le x_3<x_4$ and $x_4-x_3=x_2-x_1$. In this paper, we study the features of a more general class of Brownian motions: Brownian motions indexed by some metric space --- an $\mathbb{R}$-tree. Recall that an $\mathbb{R}$-tree is a $0$-hyperbolic metric space with desirable properties, (see \cite{Aksoy}). A detailed survey on $\mathbb R$-trees will be made in the next paragraph. Note that Brownian motion indexed by $\mathbb R$-tree is well defined. For instance, J. Istas in \cite{Istas} proved that the fractional Brownian motion (which extends Brownian motion) indexed by a hyperbolic space can be well defined  when its Hurst index $H\in(0,\frac{1}{2}]$. Furthermore, in \cite{Athreya}  the authors are able to use Dirichlet form methods to construct Brownian motion indexed by any given locally compact $\mathbb{R}$-tree. We also note that in \cite{Inoue}, it is shown that a Gaussian field (Gaussian process indexed by subset of $(\mathbb R^n,|\cdot|)$) can be represented via a set of independent increments. In this framework we study the possibility of representing a Brownian motion indexed by an $\mathbb R$-tree via the set of its independent increments. As two particular cases, we focus on $\mathbb R$-trees generated by \enquote{radial} and \enquote{river} metrics and clarify the relationship between these trees and a particular set of metric rays denoted by $\{\mathcal C_d(A,B)\}_{A,B\in M}$.  To be more precise, our investigation is motivated by the following questions:
 \begin{enumerate}
 \item Does the set of metric rays $\{\mathcal C_d(A,B)\}_{A,B\in M}$ determine the tree properties?
  \item When can an $\mathbb R$-tree be identified through the set $\{\mathcal C_d(A,B)\}_{A,B\in M}$?
  \end{enumerate}
  Since the methodology and analysis introduced in this paper are not limited to the radial and river metrics, or even to tree metrics, it is our hope that this work could lead to the interest of applying those results to Gaussian fields indexed by more general metric spaces.

The study of injective envelopes of metric spaces, also known as $\mathbb R$-trees (metric trees or $T$-theory) began with J. Tits  in \cite{Tits} in $1977$ and since then, applications have been found within many fields of mathematics.  For a complete discussion of these spaces and their relation to global  metric spaces of nonpositive curvature we refer to \cite{Bridson}.  Applications of metric trees in biology
and medicine stems from the construction of phylogenetic trees \cite{Semple}. Concepts
of \enquote{string matching} in computer science are closely related with the structure of
metric trees \cite{Bartolini}.   $\mathbb R$-trees are a generalization of an ordinary tree which allows for different
weights on edges. In order to define an $\mathbb R$-tree, we first introduce the notion of metric segment. Let $(M,d)$ be a metric space. For any $A,B\in M$, the \textit{metric segment} $[A,B]$ is defined by
$$
[A,B]=\left\{\begin{array}{ll}
\left\{X\in M:~d(A,X)+d(X,B)=d(A,B)\right\}&~\mbox{if $d(A,B)<+\infty$};\\
\emptyset&~\mbox{if $d(A,B)=+\infty$}.
\end{array}\right.
$$
In other words, $[A,B]\ne\emptyset$ if and only if $A,B$ are joined by some metric segment in $(M,d)$.
\begin{defn} [see for example \cite{Kirk}]
\label{Rtree}
An $\mathbb R$-tree is a nonempty metric space $(M,d)$
satisfying:
\begin{description}
\item[(a)] Any two points $A,B\in M$ are joined by a unique metric segment $[A,B]$.
\item[(b)] If $A,B,C\in M$, then
$
[A,B]\cap[A,C]=[A,O]~\mbox{for some}~O\in M.
$
\item[(c)]  If $A,B,C\in M$ and $[A,B]\cap[B,C]=\{B\}$, then
$
[A,B]\cup[B,C]=[A,C].
$
\end{description}
\end{defn}
There exist several different but equivalent expressions of an $\mathbb R$-tree, for more details consult \cite{Aksoy2}.  A metric space satisfying $(a)$ in Definition \ref{Rtree} is called \emph{uniquely geodesic metric space}. In the sequel we only consider uniquely geodesic metric spaces. Notice that one of the most features of an $\mathbb R$-tree is the \emph{four-point condition}. In other words, we can also characterize an $\mathbb{R}$-tree by the theorem below (see \cite{Evans}):
\begin{thm}
A uniquely geodesic metric space $(M,d)$ is an $\mathbb R$-tree if and only if it is connected, contains no triangles and satisfies the four-point condition (4PC).
\end{thm}
Recall that, $A,B,C$ form a triangle if all the triangle inequalities involving $A,B,C$ are strict and for any permutation of $(A,B,C)$, denoted by $(X,Y,Z)$, we have $[X,Y]\cap[Y,Z]=\{Y\}$. We say a metric $d$ satisfies the (4PC) if, for any $A,B, C, D$ in $M$ the following inequality holds:
 $$ d(A,B) +d(C,D) \leq \max \{d(A,C)+d(B,D), \,\,d(A,D)+d(B,C)\}.$$
The (4PC) is stronger than the triangle inequality (taking $C = D$ in the above inequality leads to the triangle inequality), but it
should not be confused with the definition of ultrametric. An ultrametric satisfies the condition
$d(A, B) \leq \max \{ d(A, C), d(B, C)\}$, and this is stronger than the (4PC).
$d$ is then said to be a tree metric if it satisfies the (4PC). Given a metric space $(M,d)$, we would capture the tree metric properties of $(M,d)$ by introducing the following sets $\{\mathcal C_d(A,B)\}_{A,B\in M}$.
\begin{defn}
\label{C}
We define, for any $P_1, P_2\in M$,
$$
\mathcal C_d(P_1,P_2)=\left\{\begin{array}{ll}
\left\{X\in M:~d(X,P_1)=d(X,P_2)+d(P_1,P_2)\right\}&~\mbox{if $d(P_1,P_2)<+\infty$};\\
\emptyset&~\mbox{if $d(P_1,P_2)=+\infty$}.
\end{array}\right.
$$
\end{defn}
Observe that two points $P_1,P_2\in M$ are joined if and only if $\mathcal C_d(P_1,P_2)\ne\emptyset$, therefore $\mathcal C_d(P_1,P_2)\ne\emptyset$ for any $P_1, P_2$ in a uniquely geodesic metric space $M$.

As one motivation, in probability theory, the sets $\{\mathcal C_d(P_1,P_2)\}_{P_1,P_2\in M}$ can be used to describe the sets of independent increments of a stochastic process. For example, let $B$ be a Brownian motion indexed by the Euclidean space $(\mathbb R^n,|\cdot|)$ in the following way:  $B(0)=0$ and the covariance structure of $B$ is given as: for $X,Y\in \mathbb R^n$,
$$
Cov(B(X),B(Y))=\frac{1}{2}\left(|X|+|Y|-|X-Y|\right).
$$
Let $d$ be the Euclidean distance defined by $d(X,Y)=|X-Y|$, then $\mathcal C_{d}(P_1,P_2)$ is precisely given by:
$$\mathcal C_{d}(P_1,P_2) = \{X\in M:\,\, B(X)-B(P_2)\,\,\, \mbox{ is independent of}\,\,\, B(P_2)-B(P_1)\}.$$
It is then of interest to ask the following questions:
\begin{description}
  \item Question $1$: Under what conditions on the set $\{\mathcal C_d(A,B)\}_{A,B\in M}$ does $(M,d)$ become an $\mathbb R$-tree?
  \item Question $2$: When can an $\mathbb R$-tree be fully identified by the set $\{\mathcal C_d(A,B)\}_{A,B\in M}$?
\end{description}
In this paper we give complete solution to Question 1 (see Section \ref{tree1} below), namely, we provide a sufficient and necessary condition on $\{\mathcal C_d(A,B)\}_{A,B\in M}$ such that $(M,d)$ is an $\mathbb R$-tree. In Section \ref{tree2}, we study Question 2 by considering radial metric and river metric. We show that the answer to Question 2 is positive for $M=\mathbb R^n$ (for some $n$) and
$$
d(A,B)=g_k(|A-B|)~\mbox{for $A,B\in \Pi_k$},
$$
where $(\Pi_k)_{k=1,\ldots,N}$ is some partition of $\mathbb R^n$ and $g_k:~\mathbb R_+\rightarrow \mathbb R_+$ is a continuous function subject to some extra properties.

\section{An Equivalence of $\mathbb R$-tree Properties}
\label{tree1}
We start by introducing the following conditions that will be used in the proof of Theorem \ref{equivtree}:\\

\textit {Condition $(A)$}:  For any 3 distinct points $A,B,C\in M$, there exists unique $O\in M$ such that
$$
\{X,Y\}\subset\mathcal C_d(Z,O)~\mbox{for any $X,Y,Z$ satisfying}~\{X,Y,Z\}\in\{A,B,C\}.
$$
Note that $(X,Y,Z)$ denotes a permutation of $(A,B,C)$.\\
\textit{Condition $(B)$}:  For any distinct $A,B,C\in M$, there exists $O\in M$ such that $$[A,B]\cap[B,C]\cap[A,C]=\{O\}.$$
Remark that if the cardinality $\# M=1$ or $2$, then $(M,d)$ is obviously an $\mathbb R$-tree, since any 2 points are joined by a unique geodesic. When $\# M\ge 3$, Condition $(A)$ guarantees that $(M,d)$ contains no circuit. The following Lemma is the key to the proof of Theorem \ref{equivtree} below:
\begin{lem}
\label{equiv}
Condition $(A)$ is equivalent to Condition $(B)$.
\end{lem}
\begin{proof}
We only consider the case where $M$ contains at least 3 distinct points. Let's pick 3 distinct points $A,B,C\in M$. Then by observing that for any distinct $X,Y\in\{A,B,C\}$,
$$
X\in\mathcal C_d(Y,O)~\mbox{is equivalent to}~O\in[X,Y].
$$
Thus Lemma \ref{equiv} holds.
\end{proof}

\begin{thm}
\label{equivtree}
A uniquely geodesic metric space $(M,d)$ is an $\mathbb R$-tree if and only if \textit{Condition $(A)$} holds.

\end{thm}
\begin{proof}
By Lemma \ref{equiv}, it is sufficient to prove that Theorem  \ref{equivtree} holds under \textit{Condition $(B)$}.
The proof consists of two steps: first we show that if $(M,d)$ is an $\mathbb R$-tree, then \textit{Condition $(B)$} is satisfied; next
 we prove that \textit{Condition $(B)$} leads to the fact that $(M,d)$ is an $\mathbb R$-tree.\\
\textbf{Step 1:} Suppose  $(M,d)$ is an $\mathbb R$-tree, since $(M,d)$ is connected, then $[A,B]\ne\emptyset$ for all $A,B\in M$. For any 3 points $A,B,C\in M$ we have:
\begin{itemize}
\item If $[A,B]\cap[B,C]=\{B\}$, then by Definition \ref{Rtree} $(c)$,
$$
\{B\}=[A,B]\cap[B,C]\subset[A,B]\cup[B,C]=[A,C].
$$
This yields
$$
[A,B]\cap[B,C]\cap[A,C]=\{B\}\cap[A,C]=\{B\}.
$$
\item If there exists $O\in M$, $O\ne B$ such that $[A,B]\cap[B,C]=[B,O]$, then $O\in[A,B]\cap[B,C]\cap[A,C]$. Thus, \textit{Condition $(B)$} is verified.
\end{itemize}
\textbf{Step 2:} Next assume \textit{Condition $(B)$} holds. By taking any $A\ne B=C$, we easily show that $[A,B]\ne\emptyset$, thus $(M,d)$ is connected. The fact that $[A,B]\cap[B,C]\cap[A,C]\ne\emptyset$ leads to the fact that there is no triangles in $(M,d)$. Then it is sufficient to prove that $d$ satisfies the (4PC). Let us pick 4 distinct points $A,B,C,D$ from $M$. Under \textit{Condition $(B)$}, there are two possibilities to the positions of $A,B,C,D$ in $(M,d)$. Namely,
\begin{enumerate}
\item Three points out of $A,B,C,D$ are in the same metric segment.
\item Case 1 above does not hold.
\end{enumerate}
 \begin{figure}[H]
\begin{minipage}[b]{0.5\linewidth}
\centering
\includegraphics[width=5cm,height=5cm]{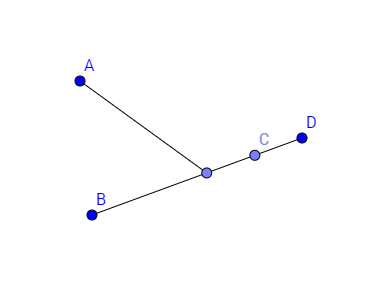}
\caption{$B,C,D$ are in one segment.}
\end{minipage}
\hspace{0.6cm}
\begin{minipage}[b]{0.5\linewidth}
\centering
\includegraphics[width=5cm,height=5cm]{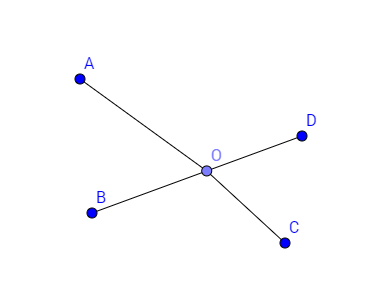}
\caption{A star graph.}
\end{minipage}
\end{figure}
In Case 1 it is easy to see that the (4PC) holds true. Indeed, without loss of generality assume $C\in[B,D]$ (see FIGURE 1), then we  necessarily have
$$
\left\{
\begin{array}{ll}
&d(A,B)\le d(A,C)+d(B,C);\\
&d(A,D)\le d(A,C)+d(D,C).
\end{array}\right.
$$
The above inequalities hold for any permutation of $A,B,C,D$. This in fact implies the (4PC). In Case 2, we observe that $A,B,C,D$ form a star graph (see FIGURE 2), i.e., there is $O\in M$ such that
 $$
 d(X,Y)=d(X,O)+d(O,Y),
 $$
 for any distinct $X,Y\in\{A,B,C,D\}$. This graph is clearly a tree hence the (4PC) is verified. Now the (4PC) is proven to be satisfied in both cases.
\end{proof}

\subsection{Characterization of $\mathcal C_d(P_1,P_2)$ for Radial Metric}
Let $(\mathbb R^n,d_1)$ $(n\ge1)$ denote an $\mathbb R$-tree with root $0$ and radial metric
$$
  d_1(A,B)=\left\{
  \begin{array}{lll}
  &|A-B|&~\mbox{if $A=aB$ for some $a\in\mathbb R$};\\
  &|A|+|B|&~\mbox{otherwise}.
  \end{array}\right.
  $$
We explicitly represent  the set $\mathcal C_{d_1}(P_1,P_2)$ for all $P_1,P_2\in\mathbb R^n$ in the following main result.
\begin{prop}
\label{Ex1}
For any $P_1,P_2\in(\mathbb R^n,d_1)$,
  \begin{equation}
  \label{Cradial}
  \mathcal C_{d_1}(P_1,P_2)=\left\{
  \begin{array}{lll}
  &[P_2,+\infty)_{\overrightarrow{0P_2}}&~\mbox{if $P_2\notin[0,P_1]_{\overrightarrow{0P_1}}$};\\
  &\mathbb R^n\backslash(P_2,+\infty)_{\overrightarrow{0P_1}}&~\mbox{if $P_2\in[0,P_1)_{\overrightarrow{0P_1}}$};\\
  &\mathbb R^n&~\mbox{if $P_1=P_2$,}
  \end{array}\right.
  \end{equation}
  where for any $A,B\in\mathbb R^n$, $[A,B)_{\overrightarrow{AB}}$ denotes the segment $\{(1-a)A+aB;a\in[0,1)\}$ and $(A,+\infty)_{\overrightarrow{0B}}$ denotes $\{aA+bB;a>1,b>0\}$ under Euclidean distance. These notations shouldn't be confused with the metric segments $[A ,B]$ of a metric space.
\end{prop}
\begin{proof}
Since it is always true that $\mathcal C_{d_1}(P_1,P_2)=\mathbb R^n$ for $P_1=P_2$, then we only consider the case when $P_1\ne P_2$. There are 3 different situations to the positions of $P_1$, $P_2$:
 \begin{description}
 \item[(1)] $P_1$, $P_2$ are on the same ray (which means, $P_2=aP_1$ for some $a\in\mathbb R$) and $0\le|P_1|<|P_2|$;
 \item[(2)]  $P_1$, $P_2$ are on the same ray and $0\le|P_2|<|P_1|$;
 \item[(3)] $P_1,P_2$ are on different rays.
\end{description}
\noindent Case $(1)$: $P_1$, $P_2$ are on the same ray and $0\le|P_1|<|P_2|$. Let $A \in C_d(P_1,P_2)$:

In this case we necessarily have
\begin{equation}
\label{d1case1}
d_1(A,P_1)=d_1(A,P_2)+|P_1-P_2|.
\end{equation}

\noindent Case $(1.1)$: If $A$ is on a different ray from $P_2$, then (\ref{d1case1}) becomes
                                                          $
                                                          |A|+|P_1|=|A|+|P_2|+|P_1-P_2|.
                                                          $
                                                          This together with the fact that $P_1\ne P_2$ implies
                                                          $
                                                          |P_1|=|P_2|+|P_1-P_2|>|P_2|.
                                                          $
                                                          This is impossible, thanks to the assumption $|P_1|<|P_2|$.

\noindent Case $(1.2)$: Suppose $A$ is on the same ray as $P_2$. Now (\ref{d1case1}) is equivalent to
                                                          $
                                                          |A-P_1|=|A-P_2|+|P_1-P_2|.
                                                          $
                                                          The solution space for $A$ is then the segment $[P_2,+\infty)_{\overrightarrow{0P_2}}$ under Euclidean distance.

                                                          We conclude that in Case $(1)$,
                                                          \begin{equation}
                                                          \label{d1case1solution}
                                                          \mathcal C_{d_1}(P_1,P_2)=[P_2,+\infty)_{\overrightarrow{0P_2}}.
                                                      \end{equation}

\noindent Case $(2)$:  $P_1,P_2$ are on the same ray and $|P_1|>|P_2|\ge0$. Note that (\ref{d1case1}) still holds.

\noindent Case $(2.1)$:  Suppose that $A$ is  on a different ray from $P_1$. (\ref{d1case1}) is then equivalent to
                                                          $
                                                          |P_1|=|P_2|+|P_1-P_2|.
                                                          $
                                                          The above equation always holds true. Therefore any $A$ on a different ray from $P_1$ belongs to $\mathcal C_{d_1}(P_1,P_2)$.

\noindent Case $(2.2)$: $A$ is on the same ray as $P_1$. Equation (\ref{d1case1}) then becomes
  $
  |A-P_1|=|A-P_2|+|P_1-P_2|,
  $
   and its solution space is segment $[0,P_2]_{\overrightarrow{0P_2}}$ under Euclidean distance.

  Combining Case (2.1) and Case (2.2), we obtain, in Case $(2)$,
  \begin{equation}
  \label{d1case2solution}
  \mathcal C_{d_1}(P_1,P_2)=\mathbb R^n\backslash(P_2,+\infty)_{\overrightarrow{0P_1}}.
  \end{equation}

                                                         \noindent Case $(3)$:  $P_1,P_2$ are on different rays, then necessarily $P_1,P_2\ne0$.

                                                          \noindent Case $(3.1)$:  $A$ is on the same ray as $P_1$.

                                                          In this case we have
                                                          $
                                                          |A-P_1|=|A|+|P_2|+|P_1|+|P_2|.
                                                          $
                                                          By the triangle inequality,
                                                          $$
                                                          |A|+|P_1|+2|P_2|=|A-P_1|\le |A|+|P_1|.
                                                          $$
                                                          This yields the absurd statement $P_2=0$!

                                                          \noindent Case $(3.2)$:  $A$ is on the same ray as $P_2$.

                                                          We have
                                                          $
                                                          |A|+|P_1|=|A-P_2|+|P_1|+|P_2|.
                                                          $
                                                          This leads to $A\in[P_2,+\infty)_{\overrightarrow{0P_2}}$.

                                                           \noindent Case $(3.3)$: $A$ is on a different ray as $P_1$, $P_2$.

                                                         In this case the fact that $
                                                          |A|+|P_1|=|A|+|P_2|+|P_1|+|P_2|$ again results a contradiction $P_2=0$.

                                                          We conclude that in Case $(3)$,
                                                          \begin{equation}
                                                          \label{d1case3solution}
                                                          \mathcal C_{d_1}(P_1,P_2)=[P_2,+\infty)_{\overrightarrow{0P_2}}.
                                                          \end{equation}

                                                        Finally, by combining (\ref{d1case1solution}), (\ref{d1case2solution}) and (\ref{d1case3solution}), we prove Proposition \ref{Ex1} holds (see FIGURES 3-4).
                                                        \end{proof}
 \begin{figure}[H]
\begin{minipage}[b]{0.45\linewidth}
\centering
\includegraphics[width=5cm,height=5cm]{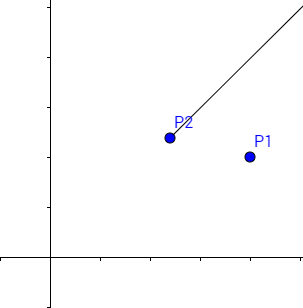}
\caption{The thick line represents the set of $\mathcal C_{d_1}(P_1,P_2)$ when $P_2$ is not in the segment $[0,P_1]$.}
\end{minipage}
\hspace{0.6cm}
\begin{minipage}[b]{0.45\linewidth}
\centering
\includegraphics[width=5cm,height=5cm]{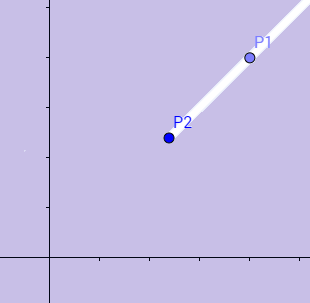}
\caption{The shaded region represents $\mathcal C_{d_1}(P_1,P_2)$ when $P_2$ is in the segment $[0,P_1)$.}
\end{minipage}
\end{figure}

                                                        Now we would show the inverse of Proposition \ref{Ex1}, namely, to answer Question 2 in Section 1: can we solve $d$ through the set of metric rays $\{\mathcal C_d(P_1,P_2)\}_{P_1,P_2\in M}$ of a given $\mathbb R$-tree $(M,d)$? For that purpose, we first state that (\ref{Cradial}) captures $\mathbb R$-tree properties.
  \begin{prop}
  \label{four}
  Let $(\mathbb R^n , d)$ be a metric space. If (\ref{Cradial}) holds for any $P_1,P_2\in(\mathbb R^n,d)$, then $(\mathbb R^n , d)$ is an $\mathbb R$-tree.
  \end{prop}
 \begin{proof} By Theorem \ref{equivtree}, we only need to show \textit{Condition $(B)$} holds.  Let us arbitrarily pick 3 different points $A,B,C\in\mathbb R^n$. If $A,B,C$ are in the same segment, saying, $A\in\mathcal C_d(B,C)$, then $C\in [A,B]\cap[B,C]\cap[A,C]$ and \textit{Condition $(B)$} is satisfied. If $A,B,C$ are not in the same segment, i.e., $X\notin\mathcal C_d(Y,Z)$ for any $\{X,Y,Z\}=\{A,B,C\}$, then we see from the definition of $\mathcal C_d(P_1,P_2)$ that
  $$
  X\in\mathcal C_d(Y,0)~\mbox{for any distinct}~X,Y\in\{A,B,C\},
  $$
  which is equivalent to $0\in[A,B]\cap[B,C]\cap[A,C]$. Hence\textit{ Condition $(B)$} is satisfied.
  \end{proof}

  \subsection{Characterization of $\mathcal C_d(P_1,P_2)$ for River Metric}
  For $A\in\mathbb R^2$, we denote by $A=(A^{(1)},A^{(2)})$. We define the $\mathbb R$-tree $(\mathbb R^2,d_2)$ with river metric $d_2$ by taking
$$
d_2(A,B)=\left\{\begin{array}{lll}
&|A^{(2)}-B^{(2)}|&~\mbox{for $A^{(1)}=B^{(1)}$;}\\
&|A^{(2)}|+|A^{(1)}-B^{(1)}|+|B^{(2)}|&~\mbox{for $A^{(1)}\ne B^{(1)}$}.
\end{array}\right.
$$
From now on we say that $A,B$ are on the same ray in $(\mathbb R^2,d_2)$ if and only if $A,B$ are on one vertical Euclidean line: $A^{(1)}=B^{(1)}$.
\begin{prop}
\label{Ex2}
Let $(\mathbb R^2,d_2)$ be a river metric space. For $P\in \mathbb R^2$, denote by $P^*=(P^{(1)},0)$ the projection of $P$ to the horizontal axis. Then for any $P_1,P_2\in\mathbb R^2$, we have
  \begin{equation}
  \label{Criver}
  \mathcal C_{d_2}(P_1,P_2)=\left\{
  \begin{array}{lll}
  &\mathbb R^2\backslash (P_2,\infty)_{\overrightarrow{P_1^*P_1}}&~\mbox{if $P_2\in[P_1^*,P_1)$};\\
  &[P_2,\infty)_{\overrightarrow{P_2^*P_2}}&~\mbox{if $P_2\notin [P_1^*,P_1)$ and $P_2^{(2)}\ne 0$};\\
  &[P_2^{(1)},\infty)_{\overrightarrow{P_1^{(1)}P_2^{(1)}}}\times\mathbb R&~\mbox{if $P_1^{(1)}\ne P_2^{(1)}$, $P_2^{(2)}=0$};\\
  &\mathbb R^2&~\mbox{if $P_1=P_2$}.
  \end{array}\right.
  \end{equation}
\end{prop}
\begin{proof}
It is obvious that $\mathcal C_{d_2}(P_1,P_2)=\mathbb R^2$ when $P_1=P_2$. For $P_1\ne P_2$, it suffices to consider 2 cases:
 \begin{description}
 \item[(1)] $P_1$, $P_2$ are on the same ray ($P_1^{(1)}=P_2^{(1)}$);
 \item[(2)] $P_1,P_2$ are on different rays ($P_1^{(1)}\ne P_2^{(1)}$).
 \end{description}
\noindent Case $(1)$:  $P_1$, $P_2$ are on the same ray. In this case we have
\begin{equation}
\label{d2case1}
d_2(A,P_1)=d_2(A,P_2)+|P_1^{(2)}-P_2^{(2)}|.
\end{equation}

 \noindent Case $(1.1)$: Suppose $A$ is on a different ray as $P_2$, then it follows from (\ref{d2case1}) that
                                                          $$
                                                          |A^{(2)}|+|P_1^{(1)}-A^{(1)}|+|P_1^{(2)}|=|A^{(2)}|+|P_2^{(1)}-A^{(1)}|+|P_2^{(2)}|+|P_1^{(2)}-P_2^{(2)}|.
                                                          $$
                                                          Since $P_1^{(1)}=P_2^{(1)}$, the above equation is simplified to
                                                          $
                                                         |P_1^{(2)}-P_2^{(2)}|+|P_2^{(2)}|-|P_1^{(2)}|=0.
                                                          $
                                                          This equation holds for all $A$ with $A^{(1)}\ne P_2^{(1)}$ provided that $P_2^{(2)}\in[0,P_1^{(2)})_{\overrightarrow{0P_1^{(2)}}}$. When  $P_2^{(2)}\notin[0,P_1^{(2)})_{\overrightarrow{0P_1^{(2)}}}$, it has no solution.

\noindent Case $(1.2)$: $A$ is on the same ray as $P_2$. Now we have
                                                          $
                                                          |A^{(2)}-P_1^{(2)}|=|A^{(2)}-P_2^{(2)}|+|P_1^{(2)}-P_2^{(2)}|.
                                                          $
                                                          The above equation holds only when $
                                                          A\in \{P_2^{(1)}\}\times[P_2^{(2)},\infty)_{\overrightarrow{0P_2^{(2)}}}$. As a conclusion, when $P_1$ and $P_2$ are on the same ray,
                                                   \begin{equation}
                                                   \label{d2case1solution}
                                                          \mathcal C_{d_2}(P_1,P_2)=\left\{\begin{array}{lll}
                                                          &\mathbb R^2\backslash (P_2,\infty)_{\overrightarrow{P_1^*P_1}}&~\mbox{if  $P_2^{(2)}\in[0,P_1^{(2)})_{\overrightarrow{0P_1^{(2)}}}$};\\
  &[P_2,\infty)_{\overrightarrow{P_2^*P_2}}&~\mbox{if  $P_2^{(2)}\notin[0,P_1^{(2)})_{\overrightarrow{0P_1^{(2)}}}$}.
                                                          \end{array}\right.
                                                          \end{equation}

                                                         \noindent  Case $(2)$: $P_1,P_2$ are on different rays.

                                                          \noindent  Case $(2.1)$: $A$ is on the same ray as $P_1$.
                                                          In this case we have
                                                          \begin{eqnarray*}
                                                          |A^{(2)}-P_1^{(2)}|&=&|A^{(2)}|+|A^{(1)}-P_2^{(1)}|+|P_2^{(2)}|+|P_1^{(2)}|+|P_1^{(1)}-P_2^{(1)}|+|P_2^{(2)}|\\
                                                          &>&|A^{(2)}|+|P_1^{(2)}|.
                                                          \end{eqnarray*}
                                                          This contradicts the triangle inequality, therefore there is no solution for $A$.

                                                          \noindent  Case $(2.2)$: $A$ is on the same ray as $P_2$.
                                                          We have
                                                          $$
                                                          |A^{(2)}|+|A^{(1)}-P_1^{(1)}|+|P_1^{(2)}|=|A^{(2)}-P_2^{(2)}|+|P_1^{(2)}|+|P_1^{(1)}-P_2^{(1)}|+|P_2^{(2)}|.
                                                          $$
                                                          By using the fact that $A^{(1)}=P_2^{(1)}$, the above equation becomes
                                                          \begin{equation}
                                                          \label{d2case22}
                                                          |A^{(2)}|=|A^{(2)}-P_2^{(2)}|+|P_2^{(2)}|.
                                                          \end{equation}
                                                          This provides:
                                                          \begin{itemize}
                                                            \item if $P_2^{(2)}=0$, then the solution space of  (\ref{d2case22}) is $ \{P_2^{(1)}\}\times \mathbb R$;
                                                            \item if $P_2^{(2)}\ne0$, then the solution space of  (\ref{d2case22}) is $
                                                            \{P_2^{(1)}\}\times[P_2^{(2)},\infty)_{\overrightarrow{0P_2^{(2)}}}$.
                                                          \end{itemize}

                                                          \noindent  Case $(2.3)$: $A$ is on a different ray from $P_1$, $P_2$. We have
                                                          \begin{eqnarray*}
                                                          &&|A^{(2)}|+|A^{(1)}-P_1^{(1)}|+|P_1^{(2)}|\\
                                                          &&=|A^{(2)}|+|A^{(1)}-P_2^{(1)}|+|P_2^{(2)}|+|P_1^{(2)}|+|P_1^{(1)}-P_2^{(1)}|+|P_2^{(2)}|.
                                                          \end{eqnarray*}
                                                          It is equivalent to
                                                          $
                                                          |A^{(1)}-P_1^{(1)}|=|A^{(1)}-P_2^{(1)}|+|P_1^{(1)}-P_2^{(1)}|+2|P_2^{(2)}|.
                                                          $
                                                          This equation has solution only when $P_2^{(2)}=0$. Provided $P_2^{(2)}=0$, the equation is written as
                                                         $
                                                          |A^{(1)}-P_1^{(1)}|=|A^{(1)}-P_2^{(1)}|+|P_1^{(1)}-P_2^{(1)}|.
                                                          $
                                                          This implies $A^{(1)}\in(P_2^{(1)},\infty)_{\overrightarrow{P_1^{(1)}P_2^{(1)}}}$.
By combining the solutions for Cases $(2.1)$, $(2.2)$, we finally obtain, in Case $(2)$,
                                                          \begin{equation}
                                                          \label{d2case2solution}
                                                          \mathcal C_{d_2}(P_1,P_2)=\left\{\begin{array}{lll}
                                                         &[P_2,\infty)_{\overrightarrow{P_2^*P_2}}&~\mbox{if $P_1^{(1)}\ne P_2^{(1)}$, $P_2^{(2)}\ne0$};\\
  &[P_2^{(1)},\infty)_{\overrightarrow{P_1^{(1)}P_2^{(1)}}}\times\mathbb R&~\mbox{if $P_1^{(1)}\ne P_2^{(1)}$, $P_2^{(2)}=0$}.
                                                          \end{array}\right.
                                                          \end{equation}
                                                       Finally, putting together Cases $(1),(2)$ completes the proof of Proposition \ref{Ex2} (see FIGURES 5-8).
                                                       \end{proof}

                                                        \begin{figure}[H]
\begin{minipage}[b]{0.45\linewidth}
\centering
\includegraphics[width=5cm,height=5cm]{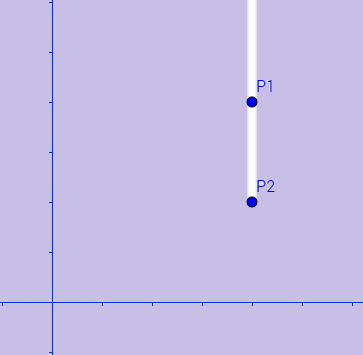}
\caption{The shaded region represents the set of $\mathcal C_{d_2}(P_1,P_2)$ when $P_2$ belongs to the segment $[P_1^*,P_1)$.}
\end{minipage}
\hspace{0.6cm}
\begin{minipage}[b]{0.45\linewidth}
\centering
\includegraphics[width=5cm,height=5cm]{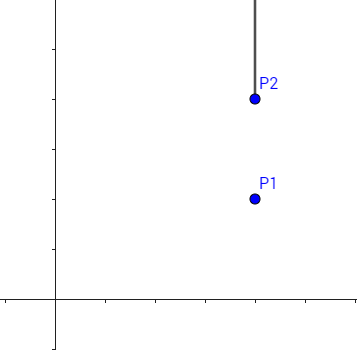}
\caption{The thick line represents $\mathcal C_{d_2}(P_1,P_2)$ when $P_1^{(1)}=P_2^{(1)}$ and $|P_2^{(2)}|>|P_1^{(2)}|$.}
\end{minipage}
\end{figure}
\begin{figure}[H]
\begin{minipage}[b]{0.45\linewidth}
\centering
\includegraphics[width=5cm,height=5cm]{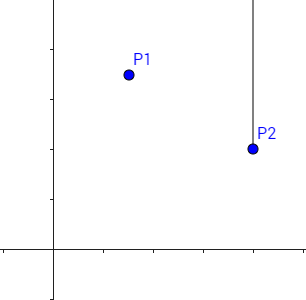}
\caption{The thick line represents the set of $\mathcal C_{d_2}(P_1,P_2)$ when $P_1^{(1)}\neq P_2^{(1)}$ and $P_2^{(2)}\ne0$.}
\end{minipage}
\hspace{0.6cm}
\begin{minipage}[b]{0.45\linewidth}
\centering
\includegraphics[width=5cm,height=5cm]{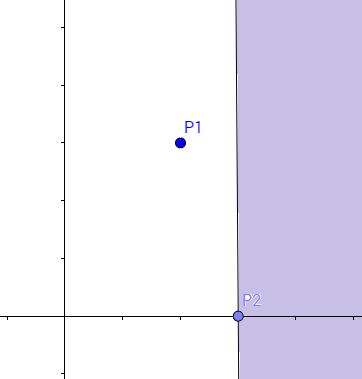}
\caption{The shaded region represents $\mathcal C_{d_2}(P_1,P_2)$ when $P_1^{(1)}\neq P_2^{(1)}$ and $P_2^{(2)}=0$.}
\end{minipage}
\end{figure}
  \begin{prop}
  \label{four1}
  Let $(\mathbb R^n,d)$ be a metric space. If
  for any $P_1,P_2\in\mathbb R^n$, (\ref{Criver}) holds, then $(\mathbb R^n,d)$ is an $\mathbb R$-tree.
  \end{prop}
 \begin{proof} We only need to show \textit{Condition $(A)$} is satisfied by the expression of $\mathcal C_d(P_1,P_2)$ in (\ref{Criver}). Observe that for any 3 distinct points $A,B,C\in\mathbb R^n$, without loss of generality, there are 3 situations according to the positions of $A,B,C$:
  \begin{description}
  \item Case $1$ : $A^{(1)}=B^{(1)}=C^{(1)}$, $A^{(2)}\in[0,B^{(2)})_{\overrightarrow{0B^{(2)}}}$, $B^{(2)}\in[0,C^{(2)})_{\overrightarrow{0B^{(2)}}}$.
      \item  Case $2$: $A^{(1)}=B^{(1)}\ne C^{(1)}$, $A^{(2)}\in[0,B^{(2)})_{\overrightarrow{0B^{(2)}}}$.
          \item Case $3$:  $A^{(1)}$, $B^{(1)}$ and $C^{(1)}$ are all distinct, $B^{(1)}\in[A^{(1)},C^{(1)}]_{\overrightarrow{A^{(1)}C^{(1)}}}$.
               \end{description}
              By (\ref{Criver}), it is easy to see \textit{Condition $(A)$} holds for $O=B$, $O=A$ and $O=(0,B^{(2)})$ respectively for Cases 1-3. Hence Proposition \ref{four1} is proven by using Theorem \ref{Rtree}.
              \end{proof}

 \section{Identification of Radial Metric and River Metric via $\mathcal C_d(P_1,P_2)$}
 \subsection{Identification of Radial Metric via $\mathcal C_{d_1}(P_1,P_2)$}
\label{tree2}
  In Proposition \ref{four} and Proposition \ref{four1}, we have shown that the sets of metric rays $\{\mathcal C_d(P_1,P_2)\}_{P_1,P_2}$ capture the tree properties of the metric spaces $(\mathbb R^n,d_1)$ and $(\mathbb R^2,d_2)$. Now we claim that subject to some additional conditions these two $\mathbb R$-trees can be uniquely identified by the sets $\{\mathcal C_d(P_1,P_2)\}_{P_1,P_2}$.
  \begin{defn}
  \label{ds}
  Let $\tilde d_1$ be a metric defined on $\mathbb R^n$ satisfying that there exists a function $f:\mathbb R_+\rightarrow\mathbb R_+$ such that
  \begin{itemize}
    \item $f$ is continuous;
    \item $f$ satisfies the following equation:

  $$
  \left\{\begin{array}{lll}
  &\tilde d_1(ax,x)=f(|ax-x|)&~\mbox{for all $x\in\mathbb R^n$ and all $a\ge0$};\\
  &f(1)=1.
  \end{array}
  \right.$$
   \end{itemize}
  \end{defn}

  \begin{thm}
  \label{inverseEx1}
   The following statements are equivalent:
  \begin{description}
  \item[(i)] $\tilde d_1=d_1$.
  \item[(ii)] For any $P_1,P_2\in(\mathbb R^n,\tilde d_1)$, $\mathcal C_{\tilde d_1}(P_1,P_2)=\mathcal C_{d_1}(P_1,P_2)$ given in (\ref{Cradial}).
  \end{description}
  \end{thm}
  Before proving Theorem \ref{inverseEx1}, we first introduce the following useful statement.
  \begin{thm} (See Aczel \cite{Aczel}, Theorem $1$)
  \label{cauchy}
  If Cauchy's functional equation
  $$
  g(u+v)=g(u)+g(v)
  $$
  is satisfied for all positive $u,v$, and if the function $g$ is
  \begin{itemize}
    \item continuous at a point;
    \item nonnegative for small positive $u-s$ or bounded in an interval,
  \end{itemize}
  then $g(u)=cu$ is the general solution for all positive $u$.
  \end{thm}
  \begin{proof} The implication $(i)\Longrightarrow(ii)$ is simply Proposition \ref{Ex1}. It remains to prove $(ii)\Longrightarrow(i)$.
  \begin{description}
  \item Case $(1)$: $A,P_1,0$ are on the same straight line with $A\ne P_1$.

 Without loss of generality, assume $|A|>|P_1|$. Then there exists $P_2\in(P_1,+\infty]_{\overrightarrow{0P_2}}$ such that $A\in[P_2,+\infty)_{\overrightarrow{0P_2}}$. By Proposition \ref{Ex1},
  $$
  \tilde d_1(A,P_1)=\tilde d_1(A,P_2)+\tilde d_1(P_1,P_2).
  $$
  Observe that $A,P_1,P_2,0$ are on the same straight line, then by the definition of $\tilde d_1$, the above equation is equivalent to
 \begin{equation}
 \label{id1}
 f(|A-P_1|)=f(|A-P_2|)+f(|P_1-P_2|).
 \end{equation}
 This is a Cauchy's equation, then by using Theorem \ref{cauchy}, the general solution is $f(u)=cu$. Together with its initial condition $f(1)=1$, we finally get $f(u)=u$. Hence,
 $$
 \tilde d_1(A,P_1)=|A-P_1|,~\mbox{for $A,P_1,0$ lying on the same straight line.}
 $$
 \item Case $(2)$: $A,P_1,0$ are not on the same straight line (in this case one necessarily has $A,P_1\ne 0$).

 We take $P_2=0$. The fact that $A\notin (0,+\infty)_{\overrightarrow{0P_1}}$ implies
 $
  \tilde d_1(A,P_1)=\tilde d_1(A,0)+\tilde d_1(P_1,0).
  $
 From Case (1) we see
 $
 \tilde d_1(A,0)=|A|~\mbox{and}~\tilde d_1(P_1,0)=|P_1|.
 $
 Therefore,
 $$
 \tilde d_1(A,P_1)=|A|+|P_1|,~\mbox{for $A,P_1,0$ lying on a different straight line.}
 $$
\end{description}
It follows from Cases (1) and (2) that $\tilde d_1=d_1$. \end{proof}
  \subsection{Identification of River Metric via $\mathcal C_{d_2}(P_1,P_2)$}

                                                        Now we claim that the inverse statement of Proposition \ref{Ex2} holds, under some extra condition.
  \begin{defn}
  \label{ds2}
   Define the metric $\tilde d_2$ on $\mathbb R^2$ by: for any $x=(x_1,x_2)$, $y=(y_1,y_2)\in\mathbb R^2$,
  $$
  \tilde d_2(x,y)=\left\{
  \begin{array}{lll}
  &g_1(|x-y|);&~\mbox{if $x_1=y_1$};\\
  &g_2(|x-y|);&~\mbox{if $x_2= y_2=0$},
  \end{array}\right.
  $$
  \end{defn}
  where $g_1,g_2$ satisfy the same conditions on $f$ in Definition \ref{ds}.
  \begin{thm}
  \label{inverseEx2}
  The following statements are equivalent:
  \begin{description}
  \item[(i)] $\tilde d_2=d_2$.
  \item[(ii)] For any $P_1,P_2\in(\mathbb R^2,\tilde d_2)$, $\mathcal C_{\tilde d_2}(P_1,P_2)=\mathcal C_{d_2}(P_1,P_2)$ given in (\ref{Criver}).
  \end{description}
  \end{thm}
  \begin{proof} The implication $(i)\Longrightarrow(ii)$ is trivial according to Proposition \ref{Ex2}. Now we prove $(ii)\Longrightarrow(i)$.
  \begin{description}
  \item Case (1):  $A^{(1)}=P_1^{(1)}$. In this case we take any $P_2\in[P_1,A]_{\overrightarrow{P_1A}}$ and get
  $
  A\in\mathcal C_{\tilde d_2}(P_1,P_2).
  $
  Equivalently,
  $$
  \tilde d_2(A,P_1)=\tilde d_2(A,P_2)+\tilde d_2(P_1,P_2).
  $$
  By using the definition of $g_1$, we obtain the following Cauchy's equation
 $$
  g_1(|A^{(2)}-P_1^{(2)}|)=g_1(|A^{(2)}-P_2^{(2)}|)+g_1(|P_1^{(2)}-P_2^{(2)}|).
  $$
Then by Theorem \ref{cauchy}, $g_1(x)=x$, for all $x\ge0$.
 \item Case $(2)$ : $A^{(1)}\ne P_1^{(1)}$.
 \item Case $(2.1)$:  We let $A^{(2)}=P_1^{(2)}=0$ and choose $P_2=(P_2^{(1)},0)$ with  $P_2^{(1)}\in[A^{(1)},P_1^{(1)}]_{\overrightarrow{P_1^{(1)}A^{(1)}}}$, then by the fact that $
 A\in\mathcal C_{\tilde d_2}(P_1,P_2)$, we have
$$
  \tilde d_2(A,P_1)=\tilde d_2(A,P_2)+\tilde d_2(P_1,P_2).
  $$
  i.e.
 $$
  g_2(|A^{(1)}-P_1^{(1)}|)=g_2(|A^{(1)}-P_2^{(1)}|)+g_1(|P_1^{(1)}-P_2^{(1)}).
  $$
 This Cauchy's equation also implies $g_2(x)=x$, for $x\ge0$.
  \item Case $(2.2)$:  $A^{(2)}\ne0$, $P_1^{(2)}=0$. In this case we take $P_2=(A^{(1)},0)$, the projection of $A$ onto the horizontal axis. Therefore by the construction of $\mathcal C_{\tilde d_2}(P_1,P_2)$ and
      \begin{eqnarray*}
      \tilde d_2(A,P_1)&=&\tilde d_2(A,P_2)+\tilde d_2(P_1,P_2)\\
      &=&|A^{(2)}|+g_2(|P_1-P_2|)\\
      &=&|A^{(2)}|+|P_1^{(1)}-A^{(1)}|.
      \end{eqnarray*}
      \item Case $(2.3)$: $A^{(2)}\ne0$, $P_1^{(2)}\ne 0$. In this case we take $P_2=(A^{(1)},0)$, the projection of $A$ onto the horizontal axis. Therefore by the construction of $\mathcal C_{\tilde d_2}(P_1,P_2)$ and
      \begin{eqnarray*}
      \tilde d_2(A,P_1)&=&\tilde d_2(A,P_2)+\tilde d_2(P_1,P_2)\\
      &=&|A^{(2)}|+|P_1^{(2)}|+|P_1^{(1)}-P_2^{(1)}|\\
      &=&|A^{(2)}|+|P_1^{(1)}-A^{(1)}|+|P_1^{(2)}|,
      \end{eqnarray*}
      we obtain that in Case $(2)$,
      $
      \tilde d_2(A,P_1)=|A^{(2)}|+|A^{(1)}-P_1^{(1)}|+|P_1^{(2)}|.
      $
      Finally for any $x,y\in\mathbb R$,
      \begin{eqnarray*}
  \tilde d_2(x,y)&=&\left\{
  \begin{array}{lll}
  &|x_2-y_2|;&~\mbox{if $x_1=y_1$};\\
  &|x_2|+|x_1-y_1|+|y_2|;&~\mbox{if $x_1\ne y_1$}
  \end{array}\right.\\
  &=&d_2(x,y).
  \end{eqnarray*}
\end{description}
\end{proof}
\section{Representation of Brownian Motion Indexed by $\mathbb R$-tree}
It should be noted that, a tree metric can be also identified by the metric segments $[A,B]$, since a uniquely geodesic metric space $(M,d)$ is a tree if and only if $[A,B]\cap[B,C]\cap[A,C]=\{O\}$ for all distinct $A,B,C\in M$. However, rather than using metric segments, the sets $\mathcal C(P_1,P_2)$ allow to capture the features of a Gaussian field, which has very important and interesting applications in the domain of random fields. As an example, Inoue and Nota (1982) \cite{Inoue} studied some classes of Gaussian fields on $(\mathbb R^n,|\cdot|)$ and represented them via the sets of independent increments. Namely, some random field $\{X(t)\}_{t\in\mathbb R^n}$ can be identified by the sets: for any $P_1,P_2\in\mathbb R^n$,
 $$
 \mathcal F_X(P_1|P_2)=\left\{A\in\mathbb R^n:~Cov\left(X(A)-X(P_2),X(P_1)-X(P_2)\right)=0\right\}.
 $$
 The set $\mathcal F_X(P_1|P_2)$ satisfies the property that,  the increments $X(A)-X(B)$ and $X(P_1)-X(P_2)$ are mutually independent if and only if $A,B\in \mathcal F_X(P_1|P_2)$. Here, we take a very similar idea of representation Gaussian fields, but  work with a tree metric which is different from Euclidean distance $|\cdot|$. More precisely,  we remark that a zero-mean Brownian motion $B$ indexed by an $\mathbb R$-tree $(M,d)$ is well-defined (see \cite{Istas,Istas2}), from its initial value $B(O)=0$ and its covariance structure
\begin{equation}
\label{Cov}
Cov(B(X),B(Y))=\frac{1}{2}\left(d(O,X)+d(O,Y)-d(X,Y)\right).
\end{equation}

Let $\{\mathcal C_d(P_1,P_2)\}_{P_1,P_2\in M}$ be the set of metric rays corresponding to $(M,d)$. Then by a similar study in \cite{Inoue}, we see that, not only $\mathcal C_d(P_1,P_2)$ can be used to identify the Brownian motion $B$, but also for any $X,Y\in \mathcal C_d(P_1,P_2)$, $B(X)-B(Y)$ and $B(P_1)-(P_2)$ are independent. This is due to the fact that, by using (\ref{Cov}) and the definition of $\mathcal C_d(P_1,P_2)$,
$$
\mathcal F_B(P_1|P_2)=\mathcal C_d(P_1,P_2),~\mbox{for any}~P_1,P_2\in M.
$$
Hence $X,Y\in\mathcal C_d(P_1,P_2)$ implies
 $
Cov\left(B(X)-B(Y),B(P_1)-B(P_2)\right)=0.
 $
 As a consequence $\{\mathcal C_d(P_1,P_2)\}_{P_1,P_2\in M}$ captures all  sets of independent increments of $\{B(X)\}_{X\in(M,d)}$. By this way one creates a new strategy to detect and simulate Brownian motion indexed by an $\mathbb R$-tree (see Section 4.2).
 \subsection{Identification of Brownian Motions Indexed by $\mathbb R$-trees}
 Let $\{B(X)\}_{X\in(\mathbb R^n,d)}$ be a zero-mean Brownian motion indexed by an $\mathbb R$-tree. Namely, $\mathbb E(B(X))=0$ for all $X\in\mathbb R^n$ and there exists an initial point $O$ such that (\ref{Cov}) holds. Then the theorems below easily follow from Theorem \ref{inverseEx1} and Theorem \ref{inverseEx2} respectively.
  \begin{thm}
  \label{inverseEx3}
   Let $B$ be a Brownian motion indexed by a metric space $(\mathbb R^n,d)$ ($n\ge1$), defined as in (\ref{Cov}). The following statements are equivalent:
  \begin{description}
  \item[(i)] $d=d_1$.
  \item[(ii)] For any $P_1,P_2\in\mathbb R^n$, $\mathcal F_B(P_1|P_2)=\mathcal C_{d_1}(P_1,P_2)$.
  \end{description}
  \end{thm}
  \begin{thm}
  \label{inverseEx4}
 Let $B$ be a Brownian motion indexed by a metric space $(\mathbb R^2,d)$, defined as in (\ref{Cov}). The following statements are equivalent:
  \begin{description}
  \item[(i)] $d=d_2$.
  \item[(ii)] For any $P_1,P_2\in\mathbb R^2$, $\mathcal F_B(P_1|P_2)=\mathcal C_{d_2}(P_1,P_2)$.
  \end{description}
  \end{thm}
\subsection{Simulation of Brownian Motion Indexed by $\mathbb R$-tree}

Let us consider a Brownian motion $B$ indexed by a tree $(\mathbb R^2,d_1)$ (recall that $d_1$ denotes radial metric) as an example. An interesting topic in statistics is to simulate such a Brownian motion. More precisely, the question is how can we generate the sample path $\{B(A_1),\ldots,B(A_n)\}$, for any different $A_1,\ldots,A_n\in(\mathbb R^2,d_1)$? In this section, we propose a new approach to simulate sample paths of Brownian motions indexed by $\mathbb R$-trees $(\mathbb R^n,d_1)$ and $(\mathbb R^2,d_2)$, which relies on the set $\mathcal C_{d_1}(P_1,P_2)$ and $\mathcal C_{d_2}(P_1,P_2)$.

The following proposition shows, in some special case, the simulation could be particularly simple.
\begin{prop}
\label{simulation}
 For any $A_1,\ldots,A_n\in(\mathbb R^2,d_1)$, there exists a permutation $\sigma\in S_n$ ($S_n$ denotes the group of permutations of $\{1,2,\ldots,n\}$) and an integer $q\ge1$ with $n_1+n_2+\ldots+n_q=n$, such that
 \begin{eqnarray}
 \label{IndepIncre}
&&\left(B(A_{\sigma(1)}),\ldots,B(A_{\sigma(n_1)})\right),\left(B(A_{\sigma(n_1+1)}),\ldots,B(A_{\sigma(n_1+n_2)})\right),\nonumber\\
&&\ldots,\big(B(A_{\sigma(n_1+\ldots+n_{q-1}+1)}),\ldots,B(A_{\sigma(n)})\big)
\end{eqnarray}
are independent, and for  each group, i.e., for $1\le l\le q$,
 \begin{equation}
 \label{sigma}
 \big(B(A_{\sigma(n_1+\ldots+n_{l-1}+1)}),\ldots,B(A_{\sigma(l)})\big)
 \end{equation}
 has independent increments.
\end{prop}
\begin{proof} It suffices to provide a such $\sigma$. We first transform $A_1,\ldots,A_n$ to their polar coordinates representations. For each $A_k$ where $k\in\{1,\ldots,n\}$, there exists $r_k\in[0,+\infty)$ and $\theta_k\in[0,2\pi)$ such that $A_k=r_ke^{i \theta_k}$. The following approach provides a permutation $\sigma$ satisfying (\ref{IndepIncre}): we choose $\sigma\in S_n$ such that
$$
\theta_{\sigma(1)}=\ldots=\theta_{\sigma(n_1)}<\theta_{\sigma_{n_1+1}}=\ldots=\theta_{\sigma_{n_1+n_2}}<\ldots<\theta_{\sigma_{n_1+\ldots+n_{q-1}+1}}=\ldots=\theta_{\sigma(n)}
$$
with $n_1+\ldots+n_q=n$ and for each group $\displaystyle \sigma\big(\sum_{m=1}^ln_m+1\big),\ldots,\displaystyle \sigma\big(\sum_{m=1}^{l+1}n_m\big)$,
$$
r_{ \sigma(\sum_{m=1}^ln_m+1)}\le r_{ \sigma(\sum_{m=1}^ln_m+2)}\le \ldots \le r_{ \sigma(\sum_{m=1}^{l+1}n_m)}.
$$
To show (\ref{IndepIncre}) and (\ref{sigma}), on one hand, by Theorem \ref{inverseEx3}, for each $l=1,\ldots,n$, the elements $\{A_k\}_{k=\sigma(\sum_{m=1}^ln_m+1),\ldots,\sigma(\sum_{m=1}^{l+1}n_m) }$ are on the same ray so they have independent increments. On the other hand, the random vectors
\begin{eqnarray*}
&&\left(B(A_{\sigma(1)}),\ldots,B(A_{\sigma(n_1)})\right),\left(B(A_{\sigma(n_1+1)}),\ldots,B(A_{\sigma(n_1+n_2)})\right),\\
&&\ldots,\big(B(A_{\sigma(n_1+\ldots+n_{q-1}+1)}),\ldots,B(A_{\sigma(n)})\big)
\end{eqnarray*}
are independent, due to the fact that for $X$, $Y$ on different rays,
\begin{eqnarray*}
Cov(B(X),B(Y))&=&\frac{1}{2}\left(d_1(X,0)+d_1(Y,0)-d_1(X,Y)\right)\\
&=&\frac{1}{2}\left(|X|+|Y|-(|X|+|Y|)\right)=0.
\end{eqnarray*}
\end{proof}
Proposition \ref{simulation} leads to the following simulation algorithm for Brownian motion indexed by $(\mathbb R^n,d_1)$.
\subsubsection{Algorithm of Simulating Brownian Motion Indexed by $(\mathbb R^2,d_1)$:}
If $A_1,\ldots,A_n$ verify the assumption given in Proposition \ref{simulation}, then
\begin{description}
\item Step $1$: Determine $\sigma\in S_n$ and $q\ge 1$ such that
\begin{eqnarray*}
&&\left(B(A_{\sigma(1)}),\ldots,B(A_{\sigma(n_1)})\right),\left(B(A_{\sigma(n_1+1)}),\ldots,B(A_{\sigma(n_1+n_2)})\right),\nonumber\\
&&\ldots,\big(B(A_{\sigma(n_1+\ldots+n_{q-1}+1)}),\ldots,B(A_{\sigma(n)})\big)
\end{eqnarray*}
are independent, and each vector has independent increments.
\item Step $2$: Generate $n$ independent zero mean Gaussian random variables $Z_1,\ldots,Z_n$, with
$$
Var(Z_k)=\left\{
\begin{array}{ll}
d_1(0,A_{\sigma(k)})&~\mbox{if $k=\sum\limits_{m=1}^{l}n_m+1$ for some $m$}\\
d_1(A_{\sigma(k-1)},A_{\sigma(k)})&~\mbox{otherwise}.
\end{array}\right.
$$
\item Step $3$: For $j=1,\ldots,n$, set
$$
B(A_{\sigma(j)})=\sum_{k=\sum_{m=1}^{l}n_m+1}^{j}Z_k,~\mbox{if $j\in \left\{\sum\limits_{m=1}^{l}n_m+1,\ldots,\sum\limits_{m=1}^{l+1}n_m\right\}$}.
$$
\end{description}
Now let us study the simulation of Brownian motion $B$ indexed by $(\mathbb R^2,d_2)$, an $\mathbb R$-tree with river metric. Similar to Proposition \ref{simulation}, we have the following proposition:
\begin{prop}
\label{simulation1}
 Given $n$ points vertically and horizontally labelled, i.e., $A_1,\ldots,A_n\in(\mathbb R^2,d_2)$ such that
 $
 \{(0,0),(A_1^{(1)},0),\ldots,(A_n^{(1)},0)\}\subset \{A_1,\ldots,A_n\}
 $
 and
\begin{eqnarray*}
&&A_{1}^{(1)}=\ldots=A_{n_1}^{(1)}<A_{n_1+1}^{(1)}=\ldots=A_{n_1+n_2}^{(1)}\\
&&<\ldots<A_{n_1+\ldots+n_{p-1}+1}^{(1)}=\ldots=A_{n_1+\ldots+n_{p}}^{(1)}< 0\\
&&\le \ldots<A_{n_1+\ldots+n_{q-1}+1}^{(1)}=\ldots=A_{n}^{(1)}
\end{eqnarray*}
with $n_1+\ldots+n_q=n$ and for $l=1,\ldots,q-1$,
$$
A_{\sum_{m=1}^ln_m+1}^{(2)}\le A_{\sum_{m=1}^ln_m+2}^{(2)}\le\ldots\le A_{\sum_{m=1}^{l+1}n_m}^{(2)}.
$$
Then there exists a sequence of independent Gaussian variables $(Z_1,\ldots,Z_{n-1})$ such that
\begin{equation}
\label{AZ}
\left(B(A_1),\ldots,B(A_n)\right)=\Big(\sum_{k\in I_1}Z_k,\ldots,\sum_{k\in I_n}Z_k\Big)~\mbox{in distribution}
\end{equation}
for some $I_k\subset\{1,\ldots,n\}$ for any $k=1,\ldots,n$.
\end{prop}
\begin{proof}  We define for $k=1,\ldots,n-1$,
\begin{equation}
\label{sigma1}
Z_k=\left\{\begin{array}{ll}
B(A_{k+1})-B(A_k)&\mbox{if $A_{k+1}^{(1)}=A_k^{(1)}$}\\
B((A_{k+1}^{(1)},0))-B((A_k^{(1)},0))&\mbox{if $A_{k+1}^{(1)}>A_k^{(1)}$}.
\end{array}\right.
\end{equation}
From  Theorem \ref{inverseEx4}, we see $(Z_1,\ldots,Z_{n-1})$ is a sequence of independent random variables. Now we are going to determine $I_1,\ldots,I_n$ such that (\ref{AZ}) holds true. Let's consider a directed graph $G=(V,E)$, with the set of vertices
$
V=\{A_1,\ldots,A_n\}
$
and the set of edges
$
E=\{e_1,\ldots,e_{n-1}\},
$
where for $k=1,\ldots,n-1$,
 \begin{equation}
\label{e}
e_k=\left\{\begin{array}{ll}
\overrightarrow{A_kA_{k+1}}&\mbox{if $A_{k+1}^{(1)}=A_k^{(1)}\ge0$}\\
\overrightarrow{(A_{k}^{(1)},0)(A_{k+1}^{(1)},0)}&\mbox{if $A_{k+1}^{(1)}>A_k^{(1)}\ge0$}\\
\overrightarrow{A_kA_{k-1}}&\mbox{if $A_{k-1}^{(1)}=A_k^{(1)}<0$}\\
\overrightarrow{(A_{k}^{(1)},0)(A_{k-1}^{(1)},0)}&\mbox{if $A_{k-1}^{(1)}<A_{k}^{(1)}<0$}.
\end{array}\right.
\end{equation}
We denote by $A_{n_0}=(0,0)$. For $k=1,\ldots,n$, let $P_k$ be the shortest path from $A_{n_0}$ to $A_k$ in $G$. Namely, there exists a set $\{k_{1},\ldots,k_{\psi(k)}\}\subset\{1,\ldots,n\}$ such that
\begin{eqnarray*}
P_{k}&=&\left(\overrightarrow{A_{n_0}A_{k_{1}}},~\overrightarrow{A_{k_{1}}A_{k_{2}}},\ldots,\overrightarrow{A_{k_{\psi(k)-1}}A_{k_{\psi(k)}}}\right)\\
&=&(e_{j_1},\ldots,e_{j_{\psi(k)}}).
\end{eqnarray*}
Denote by $I_k=\{j_1,\ldots,j_{\psi(k)}\}$, then (\ref{AZ}) is satisfied for such a choice of $(I_k)_{k=1,\ldots,n}$.
\end{proof}

From Proposition \ref{simulation1}, we provide the following simulation algorithm for Brownian motion indexed by $(\mathbb R^2,d_2)$.
\subsubsection{Algorithm of Simulating Brownian Motion Indexed by $(\mathbb R^2,d_2)$:}

If $A_1,\ldots,A_n\in \mathbb R^2$, the following algorithm shows how to simulate $(B(A_1),\ldots,B(A_n))$:
\begin{description}
\item Step $1$: Generate $n-1$ independent zero mean Gaussian random variables $Z_1,\ldots,Z_{n-1}$, with
$$
Var(Z_k)=\left\{\begin{array}{ll}
d_2(A_k,A_{k+1})&\mbox{if $A_{k+1}^{(1)}=A_k^{(1)}$}\\
d_2((A_{k+1}^{(1)},0),(A_k^{(1)},0))&\mbox{if $A_{k+1}^{(1)}>A_k^{(1)}$}.
\end{array}\right.
$$
\item Step $2$: For $k=1,\ldots,n$, determine $I_k$.
Finally,
$$
\left(B(A_k)\right)_{k=1,\ldots,n}=\Big(\sum_{j\in I_k}Z_j\Big)_{k=1,\ldots,n}~\mbox{in distribution}.
$$
\end{description}
 It is worth noting that a discrete sample path of Brownian motions indexed by $\mathbb R$-tree can be generated through its covariance matrix, where the key step is the Cholesky decomposition of the covariance matrix. Our algorithm suggests an alternative way to decompose the Brownian motion at each time step into sum of independent normal variables, with the help of $\{C_d(P_1,P_2)\}_{P_1,P_2\in\mathbb R^2}$. As an advantage to the Cholesky decomposition approach, given an $\mathbb R$-tree metric space, the sets of independent increments can be found \enquote{offline}, which will accelerate the \enquote{online} speed of our algorithms.




%
%
%
%

\noindent
\mbox{~~~~~~~}Asuman G\"{u}ven AKSOY\\
\mbox{~~~~~~~}Claremont McKenna College\\
\mbox{~~~~~~~}Department of Mathematics\\
\mbox{~~~~~~~}Claremont, CA  91711, USA \\
\mbox{~~~~~~~}E-mail: aaksoy@cmc.edu \\ \\

\noindent
\mbox{~~~~~~~}Monairah ALANSARI\\
\mbox{~~~~~~~}King Abdulaziz University\\
\mbox{~~~~~~~}Department of Mathematics\\
\mbox{~~~~~~~}Jeddah 21589, Kingdom of Saudi Arabia\\
\mbox{~~~~~~~}E-mail: malansari@kau.edu.sa \\ \\

\noindent
\mbox{~~~~~~~}Qidi PENG\\
\mbox{~~~~~~~}Claremont Graduate University\\
\mbox{~~~~~~~}Institute of Mathematical Sciences\\
\mbox{~~~~~~~}Claremont, CA  91711, USA \\
\mbox{~~~~~~~}E-mail: Qidi.Peng@cgu.edu \\ \\


\begin{thebibliography}{99}
\footnotesize
\bibitem{Aczel} {J. Aczel}, \emph{Lectures on Functional Equations and Their Applications} (Academic
Press, New York and London, 1966).
\bibitem{Aksoy}{ A. G. Aksoy and T. Oikhberg}, Some results on metric trees, Banach Center Pub. Vol.\textbf{91} (2010) 9-34.
\bibitem{Aksoy2}
     {A. G. Aksoy and S. Jin} (2014), The apple doesn't fall far from the (metric) tree: the equivalence of definitions, in \emph{Proceedings of the First Conference in Classical and Functional Analysis}, Azuga, Romania, (2014) 25-36.
\bibitem{Evans}{ S. N. Evans}, \emph{Probability and Real Trees}, Lecture Notes in Mathematics, Vol \textbf{1920 }(Springer, Berlin, 2008).
\bibitem{Athreya} {S. Athreya, M. Eckhoff and A. Winter}, Brownian motion on $\mathbb{R}$-trees. {\em Trans. of Amer. Math. Soc.}  {\bf 365, 6}  (2013) 3115-3150.
\bibitem{Bartolini} {I. Bartollini, P. Ciaccia and M. Patella}, \emph{String Matching with Metric Trees Using Approximate Distances}, Lecture Notes in Computer Science 2476 (Springer-Verlag, 2002).
 \bibitem{Bridson} {M. Bridson and A. Heefliger}, {\em Metric Spaces of Nonpositive Curvature}.  Grundlehren der Mathematischen Wissenschaften 319 ( Springer-Verlag, Berlin, 1999).
\bibitem{Inoue} {K. Inoue and A. Noda}, Independence of the increments of Gaussian random fields, {\em Nagoya Math. J.} {\bf 85} (1982) 251-268.
\bibitem{Istas} {J. Istas}, Spherical and hyperbolic fractional Brownian motion, {\em Elec. Comm. in Probabl.} {\bf 10} (2005) 254-262.
\bibitem{Istas2} {J. Istas}, Multifractional Brownian fields indexed by metric spaces with distances of negative type, {\em ESAIM.} {\bf 17} (2013) 219-223.
\bibitem{Kirk} {W.~A. Kirk }, Hyperconvexcity of $\mathbb R$-trees, {\em Fund. Math.} {\bf 156} (1998) 67-72.
\bibitem{Semple} {C. Semple and M. Steel}, \emph{Phylogenetics}, Oxford Lecture Series in Mathematics and its Applications \textbf{24} (2003).
\bibitem{Tits} {J. Tits}, A theorem  of  Lie-Kolchin for  trees,  in  \emph{Contributions  to  Algebra:  a
Collection of Papers Dedicated to Ellis Kolchin} 377-388 (Academic Press, New York, 1977).
\bibitem{Valette}  {A. Valette},  Les repr\'esentations uniformement born\'ees associ\'ees \`a un arbre, {\em Bulletin de la Soci\'et\'e Math\'ematique de Belgique S\'erie A} {\bf 42} (1990) 747-760.

\end{thebibliography}
\end{document}